\documentclass[12pt]{amsproc}
\usepackage[cp1251]{inputenc}
\usepackage[T2A]{fontenc} 
\usepackage[russian]{babel} 
\usepackage{amsbsy,amsmath,amssymb,amsthm,graphicx,color} 
\usepackage{amscd}
\def\le{\leqslant}

\newtheorem{prop}{Предложение}

\theoremstyle{definition}

\theoremstyle{plain}

\begin {document}
\centerline{УДК 512.562+517.518.85}
\unitlength=1mm
\title[Исчисление линейных расширений]
{Исчисление линейных расширений и интерполяция Ньютона}
\author{Г. Г. Ильюта}
\email{ilyuta@mccme.ru}
\address{Московский Государственный Гуманитарный Университет}
\thanks{Работа поддержана грантами РФФИ-07-01-00593 и НШ-709.2008.1.}
\maketitle

\begin{abstract}

 Мы используем разделённые разности Ньютона для вычисления сумм Грина -- рациональных функций, определяемых линейными расширениями частично упорядоченных множеств. Тождества для сумм Грина порождают соотношения для разделённых разностей Ньютона и дифференциальных форм Арнольда. Также получены обобщения интерполяционного ряда Ньютона, которые индексируются последовательностями частично упорядоченных множеств.

 We use Newton divided differences for the calculation of Greene sums -- the rational functions determined by linear extensions of partially ordered sets. Identities for Greene sums generate relations for Newton divided differences and Arnold differential forms. Also generalizations of the Newton interpolation series which are indexed by sequences of partially ordered sets are received.
\end{abstract}
 
\section {Введение}

 Сумма Грина $G(P)$ множества $P$ с отношением частичного порядка $\prec$ определяется равенством \cite{8}
$$
G(P)=\sum_{\alpha\in L(P)}\frac{1}{(x_{\alpha(1)}-x_{\alpha(2)})(x_{\alpha(2)}-x_{\alpha(3)}) \dots (x_{\alpha(n-1)}-x_{\alpha(n)})}, 
$$
где $m=|P|$ и $L(P)$ -- множество всех перестановок $\alpha\in S_m$, определяющих линейные расширения частичного порядка на $P$. Мы отождествляем элементы множества $P$ с коммутирующими переменными и рассматриваем линейное расширение как слово из этих переменных, которое начинается с некоторого максимума в $P$ и заканчивается некоторым минимумом в $P$. Подмножество частично упорядоченного множества назовём разделяющим, если его элементы попарно не сравнимы и любой элемент из дополнения либо меньше, либо больше любого элемента из этого подмножества. С помощью разделённых разностей Ньютона мы сведём задачу вычисления сумм Грина к частично упорядоченным множествам, которые могут иметь только тривиальные разделяющие подмножества -- единственный максимум, единственный минимум, и тот и другой. Для сумм Грина некоторых частично упорядоченных множеств получены явные или рекурсивные формулы. Результаты работ \cite{3}, \cite{4} и формулы этой статьи указывают на то, что со временем теория сумм Грина может стать частью теории многочленов Шуберта \cite{10}. Связи с разделёнными разностями указывают также на возможность описания сумм Грина на языке производящих функций решёточных путей и таблиц Юнга \cite{6}.

 Используя суммы Грина, мы свяжем с каждой последовательностью частично упорядоченных множеств разложение функции $(z-x)^{-1}$ в ряд из рациональных функций. Домножая на аналитическую функцию $F(z)$ и применяя теорему Коши о вычете, получим разложение в ряд функции $F(x)$. В частности, так получается интерполяционный ряд Ньютона и его предельный случай -- ряд Тейлора. Аналогично можно получить разложения для функций многих переменных. Идея построения ряда Ньютона с помощью вычета принадлежит Фробениусу \cite{7}.

 Различные разбиения множества линейных расширений $L(P)$ на подмножества приводят к различным представлениям суммы Грина $G(P)$ как суммы по подмножествам. Равенство таких сумм по подмножествам является тождеством для рациональных функций или определяемых ими функционалов, например, для разделённых разностей. Более того, подмножества в разбиении сами могут оказаться наборами линейных расширений некоторых частично упорядоченных множеств и поэтому разбиение можно рассматривать как представление частично упорядоченного множества в виде суммы других частично упорядоченных множеств. Под исчислением линейных расширений мы понимаем интерпретацию обычных тождеств как соотношений между частично упорядоченными множествами. В \cite{4} описано порождающее множество в пространстве соотношений -- эти простейшие соотношения определяются циклами в диаграммах Хассе. Ниже приводятся примеры использования этого подхода в (обобщённой) интерполяции Лагранжа. Представление частично упорядоченного множества в виде суммы цепей, отвечающих линейным расширениям, приводит к следующей формуле для интерполяционного многочлена Лагранжа $L_F(x)$ функции $F(x)$
$$
L_F(x)=f(x)\sum_{\alpha\in S_n}\frac{F(x_{\alpha(1)})}{(x-x_{\alpha(1)})(x_{\alpha(1)}-x_{\alpha(2)}) \dots (x_{\alpha(n-1)}-x_{\alpha(n)})}, 
$$
где $f(x)=\prod_{i=1}^n(x-x_i)$. Также мы покажем, что в некоторых частных случаях совпадение двух представлений мульти-функции Шура (как обобщённого определителя Якоби-Труди и как разделённой разности \cite{10}) можно интерпретировать как два способа разбить на подмножества набор линейных расширений некоторого частично упорядоченного множества.

 Тождества для сумм Грина порождают соотношения для разделённых разностей Ньютона
$$
\Delta_{X_n}[F]=\sum_{i=1}^n\frac{F(x_i)}{f'(x_i)}
$$
$$
=\frac{det(F(x_i)\quad x_i^{n-2}\quad \dots\quad x_i\quad 1)}{det(x_i^{n-1}\quad x_i^{n-2}\quad \dots\quad x_i\quad1)}\eqno (1)
$$
$$
 =\frac{1}{2\pi i}\oint\frac{F(z)}{(z-x_1) \dots (z-x_n)}dz,
$$
в определителях выписываем $i$-ю строку, $X_n=\{x_1, \dots ,x_n\}$. Cуммы Грина в тождествах можно заменить линейными комбинациями разделённых разностей функции $F$, которые получаются следующим образом: домножаем сумму Грина на $F(x_i)/\prod_{m<l<\infty}(x_i-x_l)$ и интегрируем по $x_i$. Например, так получаются классические соотношения
$$
\Delta_{X_n}[F]=\frac{\Delta_{X_n\setminus x_i}[F]-\Delta_{X_n\setminus x_j}[F]}{x_j-x_i},
$$
$$
\frac{\Delta_{X_n\setminus x_i}[F]}{(x_i-x_j)(x_i-x_k)}+\frac{\Delta_{X_n\setminus x_j}[F]}{(x_j-x_i)(x_j-x_k)}+\frac{\Delta_{X_n\setminus x_k}[F]}{(x_k-x_i)(x_k-x_j)}=0.
$$
Тождества для сумм Грина, по крайней мере по отношению к разделённым разностям, аналогичны формуле Эйлера-Якоби для вычета и детерминантному тождеству Турнбулла (более известен его частный случай -- соотношениe Грассмана-Плюккера), которые также можно использовать при доказательстве соотношений для разделённых разностей \cite{2}.

 Похожим образом тождества для сумм Грина порождают соотношения для логарифмических дифференциальных форм Арнольда
$$
\omega_{ij}=\frac{1}{2\pi i}\frac{d(x_i-x_j)}{x_i-x_j},
$$
К соотношению Арнольда \cite{1}
$$
\omega_{ij}\wedge\omega_{jk}+\omega_{ki}\wedge\omega_{ij}+\omega_{jk}\wedge\omega_{ki}=0
$$
приводит любое не линейное частично упорядоченное множество из трёх элементов. Легко видеть, что для любой перестановки $\alpha\in S_n$
$$
d(x_{\alpha(1)}-x_{\alpha(2)})\wedge d(x_{\alpha(2)}-x_{\alpha(3)})\wedge \dots \wedge d(x_{\alpha(n-1)}-x_{\alpha(n)})
$$
$$
=\sum_{i=1}^n(-1)^{n-i}dx_{\alpha(1)}\wedge \dots \wedge dx_{\alpha(i-1)}\wedge dx_{\alpha(i+1)}\wedge \dots \wedge dx_{\alpha(n)}
$$
и умножение перестановки $\alpha$ на любую транспозицию умножает это равенство на $-1$. Поэтому умножение на любую перестановку умножает равенство на знак этой перестановки и, предполагая, что исходный порядок $x_1, \dots, x_n$ является линейным расширением (единичная перестановка принадлежит $L(P)$), имеем
$$
\sum_{\alpha\in L(P)}sign(\alpha)\omega_{\alpha(1)\alpha(2)}\wedge \omega_{\alpha(2)\alpha(3)}\wedge \dots \wedge \omega_{\alpha(n-1)\alpha(n)}
$$
$$
=(\sum_{i=1}^n(-1)^{n-i}dx_1\wedge \dots \wedge dx_{i-1}\wedge dx_{i+1}\wedge \dots \wedge dx_n)G(P).
$$
По-видимому, этим фактом объясняется близость комбинаторики рациональных функций, которая появляется в \cite{11}, и комбинаторики сумм Грина. Например, формула $\tilde\Omega^{sl_2}=\Omega^{sl_2}$ на стр. 58 в \cite{11} является прямым следствием тождества Грина для частично упорядоченного множества с одним максимумом и попарно не сравнимыми остальными элементами. "Хирургия диаграмм" на стр. 65 в \cite{11} на языке статьи \cite{8} объясняется следующим замечанием: если два элемента в частично упорядоченном множестве не сравнимы, то множество линейных расширений разбивается на два подмножества согласно двум возможностям упорядочить эту пару.

\section {Обзор известных фактов о суммах Грина}

 Функция Мёбиуса $\mu:P\times P\to\mathbf Z$ частично упорядоченного множества $P$ определяется рекурсивно следующими равенствами: $\mu(a,a)=1$, $\mu(a,b)=-\sum_{a\prec c\preceq b}\mu(c,b)$, если $a\prec b$, и $\mu(a,b)=0$ в остальных случаях. Диаграммой Хассе частично упорядоченного множества называется граф, вершины которого отвечают элементам этого частично упорядоченного множества и две вершины $a$ и $b$ соединены ребром только в том случае, если $a\prec b$ и не существует вершины $c$, для которой $a\prec c\prec b$. Частично упорядоченное множество называется связным, если его диаграмма Хассе является связным графом. Частично упорядоченное множество называется плоским, если его диаграмма Хассе вкладывается с сохранением порядка в $\mathbf R\times\mathbf R$, причём, это свойство сохраняется при добавлении максимума и минимума.

 Если диаграмма Хассе частично упорядоченного множества $P$ не является связной, то $G(P)=0$ \cite{8}, \cite{3}. Для связного плоского частично упорядоченного множества $P$ имеем формулу Грина \cite{8}
$$
G(P)=\prod_{x_i\prec x_j}(x_i-x_j)^{\mu (x_i,x_j)}, 
$$
Заметим, что шаг индукции в доказательстве формулы Грина в \cite{8} сводится к следующему интерполяционному факту
$$
0=\Delta_{X_3}[x-y]
$$
$$
=\frac{x_1-y}{(x_1-x_2)(x_1-x_3)}+\frac{x_2-y}{(x_2-x_1)(x_2-x_3)}+\frac{x_3-y}{(x_3-x_1)(x_3-x_2)}, 
$$
который является частным случаем и формулы Эйлера-Якоби, и тождества Турнбулла, и самой формулы Грина. 

 Суммы Грина появились в \cite{8} как обобщения своего частного случая, связанного с диаграммами Юнга и использовавшегося при доказательстве правила Мурнагана-Накаямы для характеров симметрической группы. Модификация формулы Грина, связанная с правилом Мурнагана-Накаямы для характеров алгебры Гекке, представлена в \cite{9}.

 Если сумму Грина любого связного частично упорядоченного множества $P$ представить в виде несократимой дроби $G(P)=N(P)/D(P)$, то степень многочлена $N(P)$ равна числу независимых циклов в диаграмме Хассе (в частности, $N(P)=1$, если диаграмма Хассе является деревом), а многочлен $D(P)$ равен произведению разностей $z-x$ для элементов $x\prec z$, которые соединены ребром в диаграмме Хассе (т. е. не существует элемента $y$, для которого $x\prec y\prec z$) \cite{3}. Также в \cite{3} изучается поведение числителя $N(P)$ и знаменателя $D(P)$ при простейших преобразованиях частично упорядоченного множества, в частности, доказано, что операция стягивания ребра $(x_ix_j)$ в диаграмме Хассе следующим образом действует на сумму Грина
$$
G(P)\to\lim_{x_i\to x_j}(x_i-x_j)G(P).
$$
Любое частично упорядоченное множество можно получить с помощью стягивания рёбер в диаграмме Хассе, которая является двудольным графом (вершины можно разбить на два подмножества, в каждом из которых вершины не соединены рёбрами) \cite{3}. Тем самым, любое частично упорядоченное множество можно получить с помощью удаления и стягивания рёбер в диаграмме Хассе частично упорядоченного множества, состоящего из двух разделяющих подмножеств.

 Применение к суммам Грина интерполяционных формул стало возможным благодаря следующему свойству правой части формулы Грина и знаменателя $D(P)$ -- они представлены в виде произведений разностей переменных. В \cite{4} числитель $N(P)$ представлен в виде суммы по некоторому неоднозначно определённому набору остовных деревьев в диаграмме Хассе, причём, слагаемые в этой сумме являются произведениями разностей переменных. Некоторые результаты о разложении на множители числителя $N(P)$ и нерешённые задачи на эту тему подробно обсуждаются в \cite{4}.

 Мы уже упоминали самый существенный для понимания структуры сумм Грина результат -- описание в терминах циклов в диаграмме Хассе базисного набора в пространстве соотношений между суммами Грина \cite{4}. Интересно было бы описать пространство соотношений между соотношениями и т. д.

 Несколько примеров вычисления сумм Грина имеются в \cite{8}, \cite{3}, \cite{4}, в частности, в \cite{8} вычислены суммы Грина для стандартных диаграммам Юнга (известны функции Мёбиуса этих плоских частично упорядоченных множеств); в \cite{3} сумма Грина частично упорядоченного множества, состоящего из двух разделяющих подмножеств, представлена как мульти-функция Шура и вычислена сумма Грина в случае, когда диаграмма Хассе содержит единственный цикл; в \cite{4} получены некоторые формулы для случая двух независимых циклов в диаграмме Хассе.

\section {Суммы Грина и разделённые разности Ньютона}

 Пусть $Z_k=\{z_1,\dots,z_k\}$ является разделяющим подмножеством частично упорядоченного множества $P$ и 
$$
P^\prec=\{z\in P:z\prec z_i\},\quad P^\succ=\{z\in P:z\succ z_i\}.
$$
Для несравнимых элементов $a,b\in P$ выражение $P\cup \{a\prec b\}$ будем понимать следующим образом: к неравенствам в $P$ добавляется неравенство $a\prec b$ и все вытекающие из него по транзитивности (результат остаётся частично упорядоченным множеством). Для каждого $i\in\bar k=\{1,\dots ,k\}$ пусть
$$
P_{i\prec}=P\cup\{z_i\prec z_j:j\in \bar k\setminus i\},\quad  P_{j\succ}=P\cup\{z_j\succ z_i:i\in \bar k\setminus j\}.
$$
и для $k>1$, $i,j\in\bar k$, $i\neq j$
$$
P_{i\prec,j\succ}=P\cup\{z_i\prec z_l,z_j\succ z_l:l\in \bar k\setminus {i,j}\}.
$$
Введём также следующие обозначения
$$
R(X,Y)=\prod_{x\in X,y\in Y}(x-y),\quad p_j=\sum_{i\in\bar k}z_i^j,
$$
$$
h(z)=\prod_{i\in\bar k}(z-z_i)=R(z,Z_k),\quad Discr(h)=\prod_{i<j}(z_i-z_j)^2,
$$
$$
G_{j}^\prec=\sum_{i\in\bar k}z_i^jG(P^\prec\cup z_i),\quad G_{j}^\succ=\sum_{i\in\bar k}z_i^jG(P^\succ\cup z_i),
$$
$$
G_{j}^{\prec\succ}=\sum_{i\in\bar k}z_i^jG(P^\prec\cup z_i)G(P^\succ\cup z_i).
$$

\begin{prop}\label{prop1} Пусть $Z_k$ является разделяющим подмножеством частично упорядоченного множества $P$.

 1) Если $k=1$, то
$$
G(P)=G(P^\prec\cup Z_1)G(P^\succ\cup Z_1).
$$

 2) Если $P^\prec=\emptyset$, то
$$
G(P)=\Delta_{Z_k}[G(P^\succ\cup z)]=\sum_{j\in\bar k}\frac{G(P^\succ\cup z_j)}{h'(z_j)}.
$$

 3) Если $P^\prec=\{t\}$, то
$$
G(P)=\frac{-\Delta_{Z_k}[(t-z)^{k-1}G(P^\succ\cup z)]}{h(t)}=-\sum_{j\in\bar k}\frac{(t-z_j)^{k-1}G(P^\succ\cup z_j)}{h(t)h'(z_j)}.
$$

 4) Если $k>1$, то
$$
G(P)=-\sum_{i\neq j}\frac{(z_i-z_j)^{k-1}G(P^\prec\cup z_i)G(P^\succ\cup z_j)}{h'(z_i)h'(z_j)}
$$

 5) Если $k>3$, то
$$
G(P)=\sum_{m=0}^{k-1}(-1)^{m+1}C_{k-1}^m\Delta_{Z_k}[z^{k-m-1}G(P^\prec\cup z)]\Delta_{Z_k}[z^mG(P^\succ\cup z)]
$$
$$
=\sum_{m=0}^{k-1}\frac{(-1)^{m+1}C_{k-1}^m}{Discr(h)}\begin{pmatrix}
G_{k-1}^{\prec\succ}&G_{k+m-2}^\succ&G_{k+m-3}^\succ&\ldots&G_{m}^\succ\\
G_{2k-m-3}^\prec&p_{2k-4}&p_{2k-5}&\ldots&p_{k-2}\\
G_{2k-m-4}^\prec&p_{2k-5}&p_{2k-6}&\ldots&p_{k-3}\\
\vdots&\vdots&\vdots&\ddots&\vdots\\
G_{k-m-1}^\prec&p_{k-2}&p_{k-3}&\ldots&p_{0}
\end{pmatrix}.
$$

 При замене в 2) и 3) $P^\prec$ на $P^\succ$ и наоборот правая часть формулы умножится на $(-1)^{k-1}$ и $(-1)^k$, соответственно.
\end{prop}

 Доказательство. Из определения разделяющего множества следует, что в любом линейном расширении из $L(P)$ первыми идут все элементы множества $P^\succ$, затем все элементы множества $Z_k$ и затем все элементы множества $P^\prec$, другими словами, любая перестановка $\alpha\in L(P)$ распадается в тройку перестановок $(\alpha_\succ,\alpha_=,\alpha_\prec)$ где $\alpha_\succ\in L(P^\succ)$, $\alpha_=\in L(Z_k)(=S_k)$, $\alpha_\prec\in L(P^\prec)$ и каждая такая тройка перестановок определяет перестановку из $L(P)$. Поэтому при $k=1$ каждое слагаемое в сумме $G(P)$ является произведением слагаемого в сумме $G(P^\succ\cup Z_1)$ и слагаемого в сумме $G(P^\prec\cup Z_1)$, причём, каждая такая пара встретится. Это доказывает~1). Также имеем равенства
$$
L(P_{i\prec})=\{(\alpha_\prec,\alpha_=,\alpha_\succ)\in L(P):\alpha_=(k)=i\},
$$
$$
L(P_{j\succ})=\{(\alpha_\prec,\alpha_=,\alpha_\succ)\in L(P):\alpha_=(1)=j\},
$$
$$
L(P_{i\prec,j\succ})=\{(\alpha_\prec,\alpha_=,\alpha_\succ)\in L(P):\alpha_=(1)=j,\alpha_=(k)=i\},
$$
$$
L(P)=\cup_i L(P_{i\prec})=\cup_j L(P_{j\succ})=\cup_{i\neq j} L(P_{i\prec,j\succ}),
$$
$$
G(P)=\sum_i G(P_{i\prec})=\sum_j G(P_{j\succ})=\sum_{i\neq j} G(P_{i\prec,j\succ}).
$$

 Подмножества $\{z_i\}$, $\{z_j\}$, $\{z_i\}$ и $\{z_j\}$ являются разделяющими в частично упорядоченных множествах $P_{i\prec}$, $P_{j\succ}$, $P_{i\prec,j\succ}$, соответственно. Из 1) получаем, что в 2), 3)
$$
G(P_{j\succ})=G(P^\succ\cup z_j)G(P_{j\succ}\setminus P^\succ)
$$
и в 4)
$$
G(P_{i\prec,j\succ})=G(P^\succ\cup z_j)G(P_{i\prec,j\succ}\setminus \{P^\prec,P^\succ\})G(P^\prec\cup z_i).
$$
Частично упорядоченные множества $P_{j\succ}\setminus P^\succ$ в 2), 3) и $P_{i\prec,j\succ}\setminus \{P^\prec,P^\succ\}$ в 4) являются плоскими и по формуле Грина получаем соответствующие пунктам 2)-4) равенства
$$
G(P_{j\succ}\setminus P^\succ)=\frac{1}{h'(z_j)},
$$
$$
G(P_{j\succ}\setminus P^\succ)=-\frac{(t-z_j)^{k-1}}{h'(z_j)h(t)},
$$
$$
G(P_{i\prec,j\succ}\setminus \{P^\prec,P^\succ\})=-\frac{(x_i-x_j)^{k-1}}{h'(z_i)h'(z_j)},
$$
что доказывает 2)-4).  

 Для доказательства 5) заметим, что после сокращений в дроби $(x_i-x_j)^{k-1}/h'(z_i)h'(z_j)$ имеется множитель $(x_i-x_j)^{k-3}$. Из положительности степени следует, что в 4) можно заменить суммирование по $i\neq j$ суммированием по всем $i$ и $j$. После разложения бинома $(x_i-x_j)^{k-1}$ суммирования по $i$ и $j$ разделяются и мы получаем сумму произведений разделённых разностей. Для доказательства последнего равенства представим разделённые разности в виде отношений определителей (1), транспонируем определители в первом отношении и перемножим матрицы в числителях и матрицы в знаменателях.

 Пример 1. Если частично упорядоченное множество $P$ состоит из двух разделяющих подмножеств $Z_k$ и $X_n=\{x_1,\dots,x_n\}$, $x_j\prec z_i$, то Предложение 1 и формула Грина приводят к равенствам
$$
G(P)=(-1)^{k-1}\Delta_{Z_k}[\frac{1}{f}]=(-1)^k\Delta_{X_n}[\frac{1}{h}]
$$
$$
=(-1)^{k-1}\frac{det(\frac{1}{f(z_i)}\quad z_i^{k-2}\quad \dots\quad z_i\quad1)}{det(z_i^{k-1}\quad z_i^{k-2}\quad \dots\quad z_i\quad1)}
$$
$$
=\frac{1}{R(Z_k,X_n)}\frac{det(z_i^{k-2}f(z_i)\quad \dots\quad z_if(z_i)\quad f(z_i)\quad 1)}{det(z_i^{k-1}\quad z_i^{k-2}\quad \dots\quad z_i\quad1)},
$$
В \cite{3} выбран другой способ разбиения множества $L(P)$ на подмножества
$$
L(P)=\cup_{\alpha\in S_k}L(P\cup\{z_{\alpha(1)}\prec\dots\prec z_{\alpha(k)}\})
$$
и доказано, что $R(Z_k,X_n)G(P)$ совпадает с мульти-функцией Шура $S_{(n-1)^{k-1}}(Z_k-X_n)$. Связь мульти-функций Шура с разделёнными разностями хорошо известна \cite{10}.

 Пример 2. В случае частично упорядоченного множества, содержащего несколько разделяющих подмножеств, многократное применение Предложения 1 позволяет получить рекурсивную формулу для суммы Грина. Например, если частично упорядоченное множество $P$ состоит из трёх разделяющих подмножеств $Z_k$, $X_n$ и $Y_l=\{y_1,\dots,y_l\}$, $y_m\prec x_j\prec z_i$, то Предложение 1 и формула Грина приводят к равенствам
$$
G(P)=(-1)^{k-1}\Delta_{Z_k}[G(P^\prec\cup z)]
$$
$$
=(-1)^{k-1}\frac{det(G(P^\prec\cup z_i)\quad z_i^{k-2}\quad \dots\quad z_i\quad1)}{det(z_i^{k-1}\quad z_i^{k-2}\quad \dots\quad z_i\quad1)},
$$
где
$$
G(P^\prec\cup\{z_i\})=\frac{(-1)^n}{f(z_i)}\Delta_{Y_l}[\frac{(z_i-y)^{n-1}}{f(y)}]
$$
$$
=\frac{(-1)^n}{f(z_i)R(Y_l,X_n)}\frac{det((z_i-y_j)^{n-1}\quad y_j^{l-2}f(y_j)\quad \dots\quad y_jf(y_j)\quad f(y_j))}{det(y_j^{l-1}\quad y_j^{l-2}\quad \dots\quad y_j\quad1)}.
$$

 Пример 3. Обозначим через $P(Z_k,X_n,Y_{n-1})$ частичный порядок на множестве $Z_k\cup X_n\cup Y_{n-1}$, в котором $x_i\prec z_j$ для любых $i$, $j$ и $y_i\prec x_i$, $y_i\prec x_{i+1}$ для $i=1,\dots,n-1$. Следующий граф является диаграммой Хассе частично упорядоченного множества $P(Z_1,X_n,Y_{n-1})$

\begin{picture}(120,50)
\put(59,45){$z_1$}
\put(7,25){$x_1$}
\put(23,25){$x_2$}
\put(110,25){$x_n$}
\put(94,25){$x_{n-1}$}
\put(19,6){$y_1$}
\put(35,6){$y_2$}
\put(83,6){$y_{n-2}$}
\put(99,6){$y_{n-1}$}
\put(12,26){\line(3,1){48}}
\put(12,26){\line(1,-2){8}}
\put(28,26){\line(2,1){32}}
\put(28,26){\line(-1,-2){8}}
\put(28,26){\line(1,-2){8}}
\put(108,26){\line(-3,1){48}}
\put(108,26){\line(-1,-2){8}}
\put(92,26){\line(-2,1){32}}
\put(92,26){\line(1,-2){8}}
\put(92,26){\line(-1,-2){8}}
\put(36,10){\line(1,2){3}}
\put(84,10){\line(-1,2){3}}
\put(60,42){\line(-1,-1){5}}
\put(60,42){\line(1,-1){5}}
\put(58,23){$\dots$}
\put(60,42){\circle*{2}}
\put(12,26){\circle*{2}}
\put(28,26){\circle*{2}}
\put(92,26){\circle*{2}}
\put(108,26){\circle*{2}}
\put(20,10){\circle*{2}}
\put(36,10){\circle*{2}}
\put(84,10){\circle*{2}}
\put(100,10){\circle*{2}}
\end{picture}

\noindent Пусть
$$
g(x)=R(x,Y_{n-1}),\quad T(X_n,Y_{n-1})=\prod_{j=1}^{n-1}(x_j-y_j)(x_{j+1}-y_j).
$$
Частично упорядоченные множества $P_{i\prec}(Z_k,X_n,Y_{n-1})$ являются плоскими и, применяя Предложение 1 и формулу Грина, имеем
$$
G(P(Z_k,X_n,Y_{n-1}))=\frac{(-1)^{k-1}}{T(X_n,Y_{n-1})}\Delta_{Z_k}[\frac{g}{f}]
$$
$$
=\frac{(-1)^{k-1}}{T(X_n,Y_{n-1})}\frac{det(\frac{g(z_i)}{f(z_i)}\quad z_i^{k-2}\quad \dots\quad z_i\quad1)}{det(z_i^{k-1}\quad z_i^{k-2}\quad \dots\quad z_i\quad1)}
$$
$$
=\frac{1}{R(Z_k,X_n)T(X_n,Y_{n-1})}\frac{det(z_i^{k-2}f(z_i)\quad \dots\quad z_if(z_i)\quad f(z_i)\quad g(z_i))}{det(z_i^{k-1}\quad z_i^{k-2}\quad \dots\quad z_i\quad1)}.
$$
Формула
$$
G(P(Z_1,X_n,Y_{n-1}))=\sum_{i=1}^nG(P(Z_1,X_n,Y_{n-1})\cup\{x_i\succ x_j:j\neq i\})
$$
эквивалентна интерполяционной формуле Лагранжа (разложению на простейшие дроби)
$$
\frac{g(z_1)}{f(z_1)}=\sum_{i=1}^n\frac{g(x_i)}{(z_1-x_i)f'(x_i)}.
$$

 В следующем предложении обобщается формула из Примера~3. Результат действия симметризатора Лагранжа-Сильвестра на симметрическую функцию 
$$
F\to \sum_{I_p\subset Z_k,|I_p|=p}\frac{F(I_p)}{R(I_p,Z_k\setminus I_p)}
$$
можно представить как обобщение отношения определителей~(1) (степени переменной заменяются функциями Шура или многочленами Макдональда)~\cite{5} или как композицию элементарных разделённых разностей~\cite{10}. Предложение~1 содержит представления симметризатора Лагранжа-Сильвестра в случае, когда он действует на симметрические функции специального вида -- произведения функций от одной переменной. Пусть $Z_{i,j}=\{z_i,\dots,z_j\}$ для $i\le j$, $\Delta_{Z_{i,j}}[\frac{g}{f}]=0$ для $i>j$ и $G_{i,j}=(-1)^{j-i}G(P(Z_{i,j},X_n,Y_{n-1}))=0$ для $i>j$.

\begin{prop}\label{prop2} Для $p \le k$ справедливы равенства
$$
\frac{det(z_i^{p-1}\frac{g(z_i)}{f(z_i)}\quad \dots \quad  z_i\frac{g(z_i)}{f(z_i)}\quad\frac{g(z_i)}{f(z_i)}\quad {z_i}^{k-p-1}\quad \dots \quad z_i\quad 1)}{det(z_i^{k-1}\quad z_i^{k-2}\quad \dots\quad z_i\quad1)}\eqno (2)
$$
$$
=\sum_{I_p\subset Z_k,|I_p|=p}\frac{\prod_{z_j\in I_p}\frac{g(z_j)}{f(z_j)}}{R(I_p,Z_k\setminus I_p)}\eqno (3)
$$
$$
=det(\Delta_{Z_{i,k-p+1}}[\frac{g}{f}]\quad\Delta_{Z_{i,k-p+2}}[\frac{g}{f}]\quad\dots \quad \Delta_{Z_{i,k}}[\frac{g}{f}])\eqno (4)
$$
$$
=(T(X_n,Y_{n-1}))^p\quad det(G_{i,k-p+1}\quad G_{i,k-p+2}\quad\dots \quad G_{i,k}),\eqno (5)
$$
\end{prop}

 Доказательство. Равенство $(2)=(3)$ получается разложением определителя (2) по первым $p$ столбцам. Равенство $(3)=(4)$ доказано в \cite{12}. Равенство $(4)=(5)$ вытекает из формулы Примера~3.

 Более общие определители, включающие суммы Грина, появляются как  коэффициенты обобщённого ряда Ньютона для функции $\prod_{j=1}^p\frac{g(t_j)}{f(t_j)}$. Этот ряд можно получить как частный случай общего ряда Ньютона для функций многих переменных \cite{10}. Мы приведём прямое доказательство, используя специальный вид функции (произведение функций от одной переменной) и тождество Бине-Коши. Напомним, что в определителях выписывается $i$-я строка.

\begin{prop}\label{prop3} Разложения функции $H(z)$ в интерполяционный ряд Ньютона
$$
H(t)=\sum_{i\le j}\Delta_{X_{i,j}}[H]R(t,X_{i,j-1}),i=1,\dots ,p,
$$
влекут разложение функции $\prod_{j=1}^pH(z_j)$ в обобщённый интерполяционный ряд Ньютона
$$
\prod_{i=1}^pH(t_i)=\sum_{ j_1<\dots<j_p}det(\Delta_{Z_{i,j_1}}[H]\dots\Delta_{Z_{i,j_p}}[H])\frac{det(R(t_i,X_{j_1-1})\dots R(t_i,X_{j_p-1})}{det(1\quad t_i\quad\dots\quad t_i^{p-1})},
$$
в частности,
$$
\prod_{i=1}^p\frac{g(t_i)}{f(t_i)}=\sum_{ j_1<\dots<j_p}det(\Delta_{Z_{i,j_1}}[\frac{g}{f}]\dots\Delta_{Z_{i,j_p}}[\frac{g}{f}])\frac{det(R(t_i,X_{j_1-1})\dots R(t_i,X_{j_p-1})}{det(1\quad t_i\quad\dots\quad t_i^{p-1})}
$$
$$
=(T(X_n,Y_{n-1}))^p\sum_{ j_1<\dots<j_p}det(G_{i,j_1}\dots G_{i,j_p})\frac{det(R(t_i,X_{j_1-1})\dots R(t_i,X_{j_p-1})}{det(1\quad t_i\quad\dots\quad t_i^{p-1})}.
$$
\end{prop}

 Доказательство. Определим $(p\times\infty )$-матрицы $\Delta=(\Delta_{Z_{i,j}}[H])$, $S=(R(t_i,X_{j-1}))$ ($\Delta_{Z_{i,j}}[H])=0$ для $i>j$, $R(t_i,X_0))=1$) и применяем тождество Бине-Коши к определителю произведения $\Delta S^t$. Частный случай вытекает из формулы Примера~3. 

\section {Ряд Ньютона последовательности частично упорядоченных множеств}

 Рассмотрим последовательность связных частично упорядоченных множеств $P_1,P_2,P_3,\dots$, в каждом из которых есть элементы $x$ и $z$, которые соединены ребром в диаграмме Хассе (тогда многочлен $D(P)$ делится на $z-x$). Например, все $P_i$ могут быть одинаковыми, но их элементам соответствовуют разные переменные,

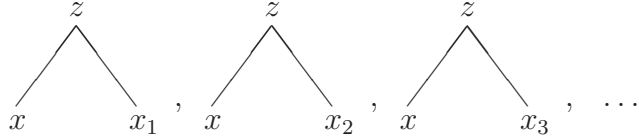
\begin{picture}(120,20)
\put(20,5){\line(3,4){8}}
\put(36,5){\line(-3,4){8}}
\put(19,2){$x$}
\put(35,2){$x_1$}
\put(27,17){$z$}
\put(41,5){$,$}
\put(46,5){\line(3,4){8}}
\put(62,5){\line(-3,4){8}}
\put(45,2){$x$}
\put(61,2){$x_2$}
\put(53,17){$z$}
\put(67,5){$,$}
\put(72,5){\line(3,4){8}}
\put(88,5){\line(-3,4){8}}
\put(71,2){$x$}
\put(87,2){$x_3$}
\put(79,17){$z$}
\put(93,5){$,$}
\put(98,5){$\dots$}
\end{picture}

 В каждой сумме Грина $G(P_i)$ разделим слагаемые на те, которые содержат множитель $(z-x)^{-1}$ и оставшиеся. Введём соответствующие обозначения
$$
G(P_i)=\frac{G'(P_i)}{z-x}=\frac{G''(P_i)}{z-x}+G'''(P_i).
$$
Тогда
$$
\frac{1}{z-x}=-\frac{G'''(P_i)}{G''(P_i)}+\frac{G'(P_i)}{G''(P_i)}\frac{1}{z-x}.\eqno (6) 
$$
Пусть
$$
G'_i=\prod_{j=1}^iG'(P_i),\quad G''_i=\prod_{j=1}^iG''(P_i).
$$
\begin{prop}\label{prop4} Последовательность связных частично упорядоченных множеств $P_1,P_2,\dots$ для всех $n>1$ определяет следующее разложение аналитической функции $F(x)$
$$
F(x)=-\sum_{i=1}^n\frac{1}{2\pi i}\oint\frac{G'_{i-1}G'''(P_i)F(z)}{G''_i}dz+\frac{1}{2\pi i}\oint\frac{G'_nF(z)}{G''_n(z-x)}dz.\eqno (7) 
$$
\end{prop}

 Доказательство. В правую часть формулы (6) для $P_1$ подставляем выражение для $(z-x)^{-1}$ из формулы (6) для $P_2$, в полученную формулу подставляем выражение для $(z-x)^{-1}$ из формулы (6) для $P_3$ и т. д. Домножая результат этих действий на $F(z)$ и интегрируя, получим формулу (7). 

 Для приведённого в первом абзаце примера последовательности трёхэлементных частично упорядоченных множеств
$$
G(P_i)=\frac{1}{(z-x)(z-x_i)},\quad G'(P_i)=\frac{1}{z-x_i},
$$
$$
G''(P_i)=\frac{1}{x-x_i},\quad G'''(P_i)=\frac{1}{(z-x_i)(x_i-x)},
$$
$$
G'_i=\prod_{j=1}^i\frac{1}{z-x_j}=\frac{1}{R(z,X_i)},\quad G''_i=\prod_{j=1}^i\frac{1}{x-x_j}=\frac{1}{R(x,X_i)},
$$
$$
-\frac{1}{2\pi i}\oint\frac{G'_{i-1}G'''(P_i)F(z)}{G''_i}dz=\Delta_{X_i}[F]R(x,X_{i-1}),
$$
$$
\frac{1}{2\pi i}\oint\frac{G'_nF(z)}{G''_n(z-x)}dz=\Delta_{X_n\cup x}[F]R(x,X_n)
$$
и формула (7) совпадает с разложением в интерполяционный ряд Ньютона
$$
F(x)=\sum_{i=1}^n\Delta_{X_i}[F]R(x,X_{i-1})+\Delta_{X_n\cup x}[F]R(x,X_n).
$$

 Рассмотрим более общий пример: частично упорядоченное множество $P_i$, $i=1,2,\dots$, состоит из элементов $z,x,x_{1i},\dots ,x_{m_ii}$ и частичный порядок определяется неравенствами
$$
x\prec z,x_{1i}\prec z,\dots ,x_{m_ii}\prec z.
$$
Пусть
$$
h_i(x)=\prod_{j=1}^{m_i}\frac{1}{x-x_{ji}}.
$$
По формуле Грина
$$
G(P_i)=\frac{h_i(z)}{z-x},\quad G'(P_i)=h_i(z),
$$
$$
G''(P_i)=h_i(x),\quad G'''(P_i)=\frac{h_i(z)-h_i(x)}{z-x},
$$
$$
G'_i=\prod_{j=1}^ih_j(z),\quad G''_i=\prod_{j=1}^ih_j(x)
$$
и формула (7) принимает вид
$$
F(x)=-\sum_{i=1}^n\frac{\frac{1}{2\pi i}\oint\prod_{j=1}^{i-1}h_j(z)\frac{h_i(z)-h_i(x)}{z-x}F(z)dz}{\prod_{j=1}^ih_j(x)}+\frac{\frac{1}{2\pi i}\oint\frac{\prod_{j=1}^nh_j(z)}{z-x}F(z)dz}{\prod_{j=1}^nh_j(x)}.
$$

 Аналогично можно обобщить интерполяционную формулу Лагранжа. Пусть в связном частично упорядоченном множестве $P$ одному из элементов соответствует переменная $z$ и несократимое представление суммы Грина имеет вид $G(P)=N(z)/D(z)$. Домножая равенство
$$
\frac{1}{z-x}=\frac{\frac{D(z)}{N(z)}-\frac{D(x)}{N(x)}}{(z-x)\frac{D(z)}{N(z)}}+\frac{\frac{D(x)}{N(x)}}{(z-x)\frac{D(z)}{N(z)}}
$$
на аналитическую функцию $F(z)$ и интегрируя по $z$, получим
\begin{prop}\label{prop5} Для связного частично упорядоченного множества имеем разложение аналитической функции $F(x)$
$$
F(x)=\frac{1}{2\pi i}\oint\frac{\frac{D(z)}{N(z)}-\frac{D(x)}{N(x)}}{(z-x)}\frac{F(z)N(z)}{D(z)}dz
$$
$$
+\frac{D(x)}{N(x)}\frac{1}{2\pi i}\oint\frac{F(z)N(z)}{(z-x)D(z)}dz,
$$
в частности, если диаграмма Хассе частично упорядоченного множества является деревом, то
$$
F(x)=\frac{1}{2\pi i}\oint\frac{D(z)-D(x)}{(z-x)}\frac{F(z)}{D(z)}dz \eqno (8)
$$
$$
+D(x)\frac{1}{2\pi i}\oint\frac{F(z)}{(z-x)D(z)}dz.
$$
\end{prop}

 Для частично упорядоченного множества $P_n$, состящего из единственного максимума  $z$ и попарно не сравнимых элементов $x_1, \dots ,x_n$, имеем равенства $N(z)=1$, $D(z)=\prod_{i=1}^n(z-x_i)$$(=f(z))$ и формула (8) представляет собой разложение функции $F(x)$ в сумму её интерполяционного многочлена Лагранжа $L_F(x)$ и остаточного члена интерполяции. Классическая формула
$$
L_F(x)=\frac{1}{2\pi   i}\oint\frac{f(z)-f(x)}{z-x}\frac{F(z)}{f(z)}dz=\sum_{i=1}^nF(x_i)\frac{f(x)}{(x-x_i)f'(x_i)} 
$$
соответствует разбиению множества $L(P_n)$ на подмножества линейных расширений, у которых вторым элементом будет $x_i$ (первым всегда будет $z$)
$$
L(P_n)=\cup_{i=1}^nL(P_n\cup\{x_j\prec x_i:j\in\bar n\setminus i\}).
$$
Разбиение множества $L(P_n)$ на одноэлементные подмножества приводит к следующей формуле

\begin{prop}\label{prop6} Для интерполяционного многочлена Лагранжа $L_F(x)$ имеем формулу
$$
L_F(x)=f(x)\sum_{\alpha\in S_n}\frac{F(x_{\alpha(1)})}{(x-x_{\alpha(1)})(x_{\alpha(1)}-x_{\alpha(2)}) \dots (x_{\alpha(n-1)}-x_{\alpha(n)})}. 
$$
\end{prop}

\bigskip


\begin{thebibliography}{74}

\bibitem{1}

В. И. Арнольд, Кольцо когомологий группы крашеных кос. Матем. заметки, 5(2), 227-231, 1969.

\bibitem{2}

Г. Г. Ильюта, Высшие порядки Брюа, формула Эйлера–Якоби и тождество Турнбулла. Успехи матем. наук, 58(4), 149-150, 2003.

\bibitem{3}

A. Boussicault, Operations on posets and rational identities of type A. Intenational Conference on Formal Power Series and Algebraic Combinatorics, 19, Tianjin, China 2007.
 
\bibitem{4}

A. Boussicault, V. Feray, Application of graph combinatorics to rational identities of type A. arXiv math. CO: 0811.2562v2.

\bibitem{5}

W. Chen, J. Louck, Interpolation for symmetric functions. Adv. Math, 117, 147–156, 1996. 

\bibitem{6}

W. Chen, B. Li, J. Louck, The flagged double Schur function. J. Alg. Comb., 15(1), 7–26, 2002.

\bibitem{7}

G. Frobenius, Uber die Entwicklung analytischer Functionen in Reihen, die nach gegebenen Functionen fortschreiten. J. reine angew. Math. 73, 1–30, 1871.

\bibitem{8}

C. Greene, A rational function identity related to the Murnaghan-Nakayama formula for the characters of $S_n$. J. Alg. Comb., 1(3), 235–255, 1992.

\bibitem{9}

T. Halverson, A. Ram, Murnaghan-Nakayama rules for characters of Iwahori-
Hecke algebras of classical type. Trans. Amer. Math. Soc., 348(10), 3967-3995, 1996.

\bibitem{10}

A. Lascoux, Symmetric Functions and Combinatorial Operators on Polynomials. CBMS Regional Conference Series in Mathematics, No. 99, AMS, Providence, 2003.
 
\bibitem{11}

R. Rimanyi, L. Stevens, A. Varchenko, Combinatorics of rational functions and Poincare-Birkhoff-Witt expansions of the canonical $U(n_{\_})$-valued differential form. Ann. Comb., 9(1), 57–74, 2005.

\bibitem{12}

H. Salzer, A determinant form for non-linear divided differences with applications. Z. Angew. Math. Mech, 66, 183–185, 1986. 

\end{thebibliography}
\end {document}